\documentclass[11pt]{article}
\usepackage{latexsym}

\newtheorem{thm}{Theorem}[section]

\newtheorem{prop}[thm]{Proposition}
\newtheorem{lem}[thm]{Lemma}

\renewcommand{\baselinestretch}{1.1} 
\def\ba{\begin{array}}
\def\ea{\end{array}}
\def\be{\begin{equation} \label}
\def\ee{\end{equation}}
\def\bea{\begin{eqnarray*}}
\def\eea{\end{eqnarray*}}
\def\beal{\begin{eqnarray} \label}
\def\eeal{\end{eqnarray}}
\def\bit{\begin{itemize}}
\def\eit{\end{itemize}}
\def\rf#1{(\ref{#1})}
\def\dot{\,\cdot\,}
\def\proof{{\em Proof. }}

\def\a{\alpha}
\def\g{\gamma}

\def\d{\delta}
\def\D{\Delta}
\def\e{\varepsilon}
\def\l{\lambda}
\def\L{\Lambda}
\def\r{\varrho}
\def\s{\sigma}

\def\Th{\Theta}

\def\O{\Omega}

\def\LL{{\cal L}}

\def\G{{\cal G}}

\def\Z{{\bf Z}}
\def\N{{\bf N}}
\def\R{{\bf R}}
\def\ti{\to\infty}
\def\c{\circ}
\def\b{\bullet}

\def\st{{\mbox{\tiny\rm stag}}}
\def\even{{\mbox{\tiny\rm even}}}
\def\odd{{\mbox{\tiny\rm odd}}}
\def\ord{{\mbox{\tiny\rm ord}}}
\def\per{{\mbox{\tiny\rm per}}}
\def\dis{{\mbox{\tiny\rm dis}}}
\def\sea{S}
\def\inv{^{-1}}

\begin{document}

\title{{\bf Entropy-driven phase transitions\\ in multitype lattice 
gas models}}
\author{%
Hans-Otto Georgii\\ {\small\sl Mathematisches Institut der 
Universit\"at M\"unchen}\\ {\small\sl Theresienstr.\ 39, D-80333 
M\"unchen, Germany}
\and Valentin Zagrebnov\\{\small\sl Universit\'e de la M\'editerran\'ee 
and  Centre de Physique Th\'eorique,}\\
{\small\sl CNRS-Luminy-Case 907, 13288 Marseille Cedex 9, France}} 

\date{}
\maketitle

\renewcommand{\baselinestretch}{1.0}
\small
\begin{quote}
\date{February 23, 2000}    \\[-0.6ex]
\hspace*{0ex}\hrulefill\hspace*{0ex}

In multitype lattice gas models with hard-core interaction of
Widom--Rowlinson type, there is a competition between the entropy due
to the large number of types, and the positional energy and geometry
resulting from the exclusion rule and the activity of particles.  We
investigate this phenomenon in four different models on the square
lattice: the multitype Widom--Rowlinson model with diamond-shaped
resp.\ square-shaped exclusion between unlike particles, a
Widom--Rowlinson model with additional molecular exclusion, and a
continuous-spin Widom--Rowlinson model.  In each case we show that
this competition leads to a first-order phase transition at some
critical value of the activity, but the number and character of phases
depend on the geometry of the model.  Our technique is based on
reflection positivity and the chessboard estimate.
\\[-0.6ex]
\hspace*{0ex}\hrulefill\hspace*{0ex}

\end{quote}
\renewcommand{\baselinestretch}{1.1}\normalsize

\section{Introduction}

Although the most familiar examples of phase transitions in lattice
models originate from a degeneracy of ground states and therefore
occur at low temperatures, this is not the only situation in which
phase transitions can occur.  Another possible source of criticality
is a conflict of energy and entropy.  This was noticed first by
Dobrushin and Shlosman \cite{DS} in the case of an asymmetric
double-well potential with two (sharp resp.\ mild) local minima
separated by a barrier.  They found that for some specific temperature
energy and entropy attain a balance leading to the coexistence of
high- and low-temperature phases corresponding to the two wells; cf.\
also Section 19.3.1 of \cite{Gii}.  Later on, Koteck\'y and Shlosman
\cite{KS} observed that such a first-order phase transition can occur
even in the absence of an energy barrier, provided there is an
``explosion'' of entropy.  They demonstrated this in particular on the
prototypical case of the $q$-state Potts model on $\Z^d$ for large
$q$, showing that for a critical temperature there exist $q$
distinct ordered low-temperature phases as well as one disordered
high-temperature phase; see also Section 19.3.2 of \cite{Gii}.

This paper has the objective of studying entropy-driven first-order phase
transitions of similar kind in multitype lattice gas models with
type-dependent hard core interaction.  In such models the crucial
parameter is the activity instead of temperature, and the
entropy-energy conflict turns into a competition between the
entropy of particle types and the positional energy and geometry
resulting from the exclusion rule and the activity of particles.  One
is asking for a critical activity with coexistence of low-density and
high-density phases.

The basic example of this kind is the multicomponent Widom--Rowlinson
lattice gas model investigated first by Runnels and Lebowitz \cite{RL}
in 1974 and studied later (theoretically and numerically) by Lebowitz
et al.\ \cite{LMNS}, cf.\ also \cite{NL}.  If the number $q$ of types
is large enough (the numerical estimates give $q\ge 7$), there exist
three different regimes: besides the low-density uniqueness regime and
a high-density regime with $q$ ``demixed'' phases for $z>z_c(q)$,
there exists an intermediate domain of activities $z_0(q)<z<z_c(q)$ with
two ``crystal'' (or ``staggered'') phases with an occupation pattern of 
chessboard type, and the phase transition at
$z_c(q)$ is of first order.   
The transition between staggered and demixed phases at
$z_c(q)$ is again entropy-driven: in the staggered phases the type
entropy wins, with the effect of an entropic repulsion of positions
forcing the particles onto a sublattice, whereas in the demixed
phases the particles gain energy and positional freedom but loose their type
entropy.  The same kind of phenomenon has also been discovered for a
class of spin systems with annealed dilution (including diluted Potts
and plane rotor models) \cite{CKS,CKS2}.

The aim of the present paper is to analyze the interplay of type 
entropy and the geometry induced by the lattice and the exclusion rule.
While we stick to the integer lattice $\Z^d$ (and for simplicity in 
fact to the case $d=2$), we vary the exclusion rule in order to gain 
some insight into the geometric effects involved. We investigate and 
compare four different models:
\begin{enumerate}
    \item the standard multitype Widom--Rowlinson model; 
    \item a multitype Widom--Rowlinson model with nearest-neighbor and
    next-nearest neighbor exclusion between particles of different type;
    \item a multitype Widom--Rowlinson model with additional
    type-independent hard-core interaction between next-nearest particles;
    \item a ferrofluid model of oriented particles with exclusion between
    neighboring but not sufficiently aligned particles.  (Similar
    continuous-spin counterparts of models 2 and 3 will also be
    considered.)
\end{enumerate}
We show that in each of these examples an entropy-driven first-order
phase transition occurs, but the number and specific characteristics
of coexisting phases are different in all cases.  Our technique is
similar to that used in \cite{CKS,CKS2,DS,KS} and Chapters 18/19 of
\cite{Gii}, and is based on the (trivial) reflection positivity in
lines through lattice sites and the resulting chessboard estimate
\cite{Gii}.  While in model 1 this is only an alternative (and perhaps
more elementary) approach to the results obtained in \cite{LMNS} by
means of Pirogov-Sinai theory, the very same argument works also in
the other models with only slight modifications.

This paper is organized as follows.  In Section \ref{results} we
introduce the four models and present our results.  The proofs follow
in Sections \ref{sec:pf_WR} to \ref{sec:pf_cs}.  The general scheme is
explained in detail for model 1, the standard Widom--Rowlinson lattice
model.  In the other cases we only indicate the necessary
changes.
\medskip

\renewcommand{\baselinestretch}{1.0}
\small
{\em Acknowledgement}. H.O.G. gratefully acknowledges warm hospitality of 
the Centre de Physique Th\'eorique in Marseille-Luminy, and V.Z. of the 
Mathematical Institute of the University of Munich. This work was 
supported by the Deutsche Forschungsgemeinschaft, SPP 1033.
\renewcommand{\baselinestretch}{1.1}
\normalsize

\section{Models and results}
\label{results}

\subsection{The multitype Widom--Rowlinson lattice gas}
\label{sec:WR}

This model describes a system of particles of $q$ different types
(`colors') which are allowed to sit on the sites of the square lattice
$\Z^2$.  (For simplicity we stick to the two-dimensional case; an
extension to higher dimensions is straightforward, cf.\ Chapter 18 of
\cite{Gii} or \cite{Fuku}.)  At each lattice site we have a random 
variable $\s_i$
taking values in the set $E=\{0,1,\ldots,q\}$.  The equality $\s_i=0$
means that site $i$ is empty, and $\s_i=a\in\{1,\ldots,q\}$ says that
$i$ is occupied by a particle of color $a$.  Particles of different
color interact by a hard-core repulsion: they are not allowed to sit
next to each other. There is no interaction between particles of the 
same color. This means that the formal Hamiltonian has the
form
\be{Ham}
H(\s)= \sum_{\langle ij\rangle} U(\s_i,\s_j)\;,
\ee
where the sum extends over all nearest-neighbor pairs $\langle 
ij\rangle\subset\Z^2$ of lattice sites (i.e., $|i-j|=1$), and the 
potential $U$ is given by
\be{U}
U(\s_i,\s_j)=\left\{\ba{cl}\infty&\mbox{if }0\ne \s_i\ne \s_j\ne 0 
\;,\\
0&\mbox{otherwise.}\ea
\right.
\ee
This model is a lattice analog of the continuum two-species model of 
Widom and Rowlinson \cite{WR}, and was first introduced by Lebowitz 
and Gallavotti \cite{LG} (for $q=2$) and Runnels and Lebowitz 
\cite{RL} (for general $q$).

Since $U$ is either $0$ or $\infty$, the temperature does not play 
any role, and the only energetic parameter is the activity $z>0$ 
which governs the overall particle density; we assume that the activity 
does {\em not\/} depend on the particle color. Accordingly, the Gibbs 
distribution in a finite region $\L\subset\Z^2$ with boundary 
condition $\eta$ in $\L^c=\Z^2\setminus\L$ is given by
\be{Gibbs}
\mu_{\L,\eta}^{z,q}(\s)= 1_{\{\s\equiv\eta \mbox{ \scriptsize off }\L\}}\;
(Z_{\L,\eta}^{z,q})^{-1}\; z^{N_\L(\s)} \exp\Big[-
\sum_{\langle ij\rangle\cap\L\ne\emptyset} U(\s_i,\s_j)\Big]\;,
\ee
where $N_\L(\s)=|\{i\in\L:\s_i\ne 0\}|$ is the number of particles 
in $\L$, and $Z_{\L,\eta}^{z,q}$ is a normalizing constant.

Alternatively, we may think of $\mu_{\L,\eta}^{z,q}$ as obtained by 
conditioning a Bernoulli measure on the set of admissible 
configurations. Let
\[
\O=\{ \s\in E^{\Z^2}: \forall \langle ij\rangle\ 
\s_i\s_j=0\mbox{ or }\s_i=\s_j\}
\]
be the set of all admissible configurations on $\Z^2$. Given any
such configuration $\s\in\O$ and any subset $\L$ of $\Z^2$, we write
$\s_\L$ for the restriction of $\s$ to $\L$. We also write
$\O_{\L,\eta}$ for the set of all 
admissible configurations $\s\in E^{\L}$ in $\L$ which are compatible 
with some $\eta\in\O$, in the sense that the composed 
configuration $\s\eta_{\L^c}$ belongs to $\O$. In particular, we write
$\O_\L=\O_{\L,0}$ for the set of all admissible configurations in 
$\L$, which are compatible with the empty configuration $0$ outside $\L$.
It is then easy to see that
\[
\mu_{\L,\eta}^{z,q} = \pi_{\L}^{z,q}(\dot|\O_{\L,\eta})\;,
\]
where $\pi_{\L}^{z,q}=\bigotimes_{i\in\L} \pi_i^{z,q}$ is the $\L$-product
of the measures 
$\pi_i^{z,q}=(\frac1{1+qz},\frac{z}{1+qz},\ldots,\frac{z}{1+qz})$ on 
$E$. 

Given the Gibbs distributions $\mu_{\L,\eta}^{z,q}$, we define the 
associated class $\G(z,q)$ of (infinite volume) Gibbs measures on $\O$ in the 
usual way \cite{Gii}. Our main result below shows that for large $q$
there exist two different activity regimes in which $\G(z,q)$ 
contains several phases of quite different behavior. These regimes 
meet at a critical activity $z_c(q)$ and produce a first-order phase 
transition. 

The different phases admit a geometric description in percolation
terms.  Let $\Z^2$ be equipped with the usual graph structure
(obtained by drawing edges between sites of Euclidean distance 1). 
Given any $\s\in\O$, a subset $S$ of $\Z^2$ will be called an
\emph{occupied cluster} if $S$ is a maximal connected subset of
$\{i\in\Z^2: \s_i\ne 0\}$, and an \emph{occupied sea} if, in addition,
each finite subset $\D$ of $\Z^2$ is surrounded by a circuit (i.e.,
closed lattice path) in $S$.  In other words, an occupied sea is an infinite
occupied cluster with interspersed finite `islands'.  If in fact
$\s_i=a$ for all $i\in S$ we say $S$ is an \emph{occupied sea of color
$a$}.  We consider also the dual graph structure of $\Z^2$ with
so-called $*$edges between sites of distance 1 or $\sqrt 2$, and the
associated concept of $*$connectedness.  An \emph{even occupied
$*$sea} is a maximal $*$connected subset of $\{i=(i_1,i_2)\in\Z^2:
i_1+i_2 \mbox{ even, } \s_i\ne 0\}$ containing $*$circuits around
arbitrary finite sets $\D$.  Likewise, an \emph{odd empty $*$sea} is a
maximal $*$connected subset of $\{i=(i_1,i_2)\in\Z^2: i_1+i_2 \mbox{
odd, } \s_i=0\}$ surrounding any finite $\D$.

\begin{thm}\label{th:WR}
If the number $q$ of colors exceeds some $q_0$, there exists an
activity threshold $z_c(q)\in \;]\,q/5,5q\,[$ and numbers
$0<\e(q)<1/3$ with $\e(q)\to 0$ as $q\ti$ such that the following
hold:

{\rm (i) }For $z>z_c(q)$, there exist $q$ distinct translation
invariant `colored' phases $\mu_a\in\G(z,q)$, $a\in\{1,\ldots,q\}$. 
Relative to $\mu_a$, there exists almost surely an occupied sea of
color $a$ containing any given site with probability at least
$1-\e(q)$.

{\rm (ii) }For $q_0/q\le z<z_c(q)$, there exist two distinct
`staggered' phases $\mu_\even,\mu_\odd\in\G(z,q)$ invariant
under even translations.  Relative to $\mu_\even$, there exist almost
surely both an even occupied $*$sea and an odd empty $*$sea, and any
two adjacent sites belong to these $*$seas with probability at least
$1-\e(q)$.  In addition, all occupied clusters are finite almost
surely, and their colors are independent and uniformly distributed
conditionally on their position.  $\mu_\odd$ is obtained from
$\mu_\even$ by a one-step translation.

{\rm (iii) }At $z=z_c(q)$, a first-order phase transition occurs, in 
the sense that $q+2$ distinct phases $\mu_\even,\mu_\odd,\mu_1,\ldots, 
\mu_q\in\G(z_c(q),q)$ coexist which enjoy the properties above.
\end{thm}

The preceding theorem can be summarized by the following phase diagram.

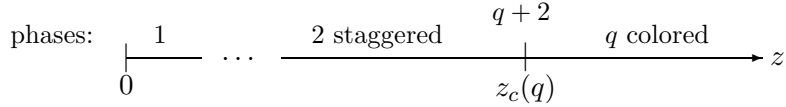
\begin{figure}[htb]
\begin{center}
\setlength{\unitlength}{1cm}
\begin{picture}(7.6,1)
\put(0.14,0){\line(1,0){1}}
\put(1.4,-0.024){$\ldots$}
\put(2.2,0){\vector(1,0){6.4}}
\put(8.7,-0.1){$z$}
\put(0,-0.31){$\mbox{$|$}\atop\mbox{0}$}
\put(5,-0.32){$\mbox{$|$}\atop\mbox{$z_c(q)$}$}
\put(-1.4,0.2){\small phases:}
\put(0.5,0.2){\small 1}
\put(2.6,0.2){\small 2 staggered}
\put(6.5,0.2){\small $q$ colored}
\put(5,0.5){\small $q+2$}
\end{picture}
\end{center}
\end{figure}
We continue with a series of comments.

\medskip\noindent
{\bf Remark \thethm } (1) The existence of staggered phases in an
intermediate activity region was first observed by Runnels and
Lebowitz \cite{RL}.  As will become apparent later, this is
a consequence of the fact that the lattice $\Z^2$ is bipartite and the
interaction is nearest-neighbor.  According to Theorem \ref{th:WR},
for large $q$ the staggered regime extends up to the fully ordered
regime, and the transition from the staggered regime to the ordered
regime at $z_c(q)$ is of first order.  This result (which disproves a
conjecture in \cite{RL}) has already been obtained before by Lebowitz,
Mazel, Nielaba and \v{S}amaj \cite{LMNS}.  While their argument relies on
Pirogov--Sinai theory (which even gives the asymptotics of $z_c(q)$), 
we offer here a different proof based on
reflection positivity which is quite elementary and
can easily be adapted to our other models
(including a continuous-spin variant of the present model).

\smallskip\noindent 
(2) There are two kinds of ordering to be distinguished: positional
order and color-order.  The colored (or `demixed') phases
$\mu_1,\ldots, \mu_q$ show color-order but no positional order. (The 
impression of positional order is a delusion coming from the lattice 
regularity.) Their high density takes advantage of the chemical 
energy of particles (i.e., of the activity $z$.) On the other hand, the
staggered (or `crystal') phases $\mu_\even,\mu_\odd$ exhibit
positional order but color-disorder.  Positional \emph{and} color-disorder 
occurs in the uniqueness regime at sufficiently low activities.

\smallskip\noindent 
(3) The first-order transition at $z_c(q)$ manifests itself
thermodynamically by a jump of the particle density as a function of
the activity.  In fact, $z_c(q)$ can be
characterized as the unique value where the density jumps over the
level $2/3$, cf.\ Lemma \ref{z_c}.

\smallskip\noindent
(4) For small $z$ there exists only one Gibbs measure in $\G(z,q)$. 
For example, using disagreement percolation one easily finds that this
is the case when $qz<p_c/(1-p_c)$, where $p_c$ is the Bernoulli site
percolation threshold for $\Z^2$; see \cite{vdBM,GHM} for more
details.  We do not know whether the uniqueness regime extends right up to
the staggered regime.  As will be explained in the next comment, this
question is related to the behavior of the hard-core lattice gas.

\smallskip\noindent
(5) As was already noticed in \cite{RL,LMNS}, the occupation structure
of the large-$q$ Widom--Rowlinson model at activity $z$ is
approximately described by the hard-core lattice gas with activity
$\zeta=qz$.  This becomes evident from the following argument (which
is more explicit than those in \cite{RL,LMNS}).  Consider the Gibbs
distribution $\mu_{\L,\per}$ of our model in a rectangular box $\L$
with periodic boundary condition.  (This boundary condition is
natural, since later on we will only look at phases appearing in the
extreme decomposition of infinite volume limits of $\mu_{\L,\per}$; we
could also use empty or monochromatic boundary conditions instead.)  Let
$\psi_{\L,\per}$ be the image of $\mu_{\L,\per}$ under the projection
$\s\to n=(n_i)_{i\in\L}=(1_{\{\s_i\ne 0\}})_{i\in\L}$ from $E^\L$ to
$\{0,1\}^\L$ mapping a configuration of colored particles onto the
occupation pattern.  $\psi_{\L,\per}$ is called the
site-random-cluster distribution, see Section 6.7 of \cite{GHM}.  Its
conditional probabilities are given by the formula
\[
\psi_{\L,\per}(n_i=1|n_{\L\setminus\{i\}}) = \frac{qz}{qz + q^{\kappa(i,n)}}\;,
\]
where $\kappa(i,n)$ is the number of clusters of
$\{j\in\L\setminus\{i\}:n_j=1\}$ meeting a neighbor of $i$.  Now,
since $\kappa(i,n)=0$ if and only if all neighbors of $i$ are empty,
this conditional probability tends to the one of the hard-core lattice
model with activity $\zeta$ when $q\ti$ and $\zeta=qz$ stays fixed. 
Unfortunately, this result is \emph{not} sufficient to conclude that in this
limit the transition point from the unique to the staggered phase
converges to that of the hard-core lattice gas, although this seems 
likely and is suggested by simulations \cite{LMNS}.
 
\smallskip\noindent
(6) One may ask whether the monotonicity of the transition from the
staggered to the ordered regime can be deduced from stochastic
monotonicity properties of the site-random-cluster model, as is
possible in the Potts model.  Unlike in the standard (bond)
random-cluster model, such a stochastic monotonicity is not available
\cite{BHW}.  To obtain the existence of a unique transition point
$z_c(q)$ we will therefore use the convexity of the pressure, which
implies that the particle density is an increasing function of $z$.

\smallskip\noindent
(7) As is often the case in this kind of context, our bounds on $q_0$
are not very useful.  They only allow us to conclude that we can take
e.g.\ $q_0= 2\cdot 10^{85}$.  Note that for small $q$ the ordered regime
still exists \cite{RL}, but for $q=2$ there is no staggered regime but
instead a direct second order transition from the gas phase to the
ordered phase.  For this and more information about the minimal $q$ at
which the staggered phase appears see \cite{LMNS,NL}.

\subsection{The square-shaped Widom--Rowlinson lattice gas}
\label{sec:sqWR}

The standard Widom--Rowlinson model considered above is defined by the 
exclusion rule that no two particles of different color may occupy 
adjacent sites. Equivalently, one may think of the particles as having 
the shape of the diamond $\{x\in\R^2: \|x\|_1<1\}$, and diamonds of 
different color are required to be disjoint.

In this section we want to study a variant with different geometry: we
identify a particle at position $i\in\Z^2$ with the suitably colored
square $\{x\in\R^2: \|x-i\|_\infty<1\}$, and we stipulate that
squares of different color are disjoint, while squares of the same 
color may overlap. Alternatively, this 
assumption amounts to replacing $\Z^2$ by its matching dual 
$(\Z^2)^*$, which is obtained from $\Z^2$ by keeping all 
nearest-neighbor bonds and adding bonds 
between diagonal neighbors of Euclidean distance $\sqrt 2$; for lack 
of a generally accepted name we call this lattice the \emph{face-crossed 
square lattice}. Accordingly, we say that two sites $i,j\in\Z^2$ 
are $*$adjacent if $|i-j|=1$ or $\sqrt 2$, and we write $\langle 
ij\rangle^*$ for such pairs of sites. 
Saying that squares of different colors are disjoint 
is then equivalent to saying that particles of different color do
not sit on $*$adjacent sites. 

The point of considering this model is that the square-shape of 
particles fits better with the geometry of $\Z^2$ than the 
diamond-shape in the standard Widom--Rowlinson model. As a 
consequence, the somehow artificial staggered phases disappear, 
and the model is closer to what one would expect to hold in the 
continuum Widom--Rowlinson model.

We now turn to precise statements. The formal Hamiltonian of the 
square-shaped Widom--Rowlinson model is analogous to \rf{Ham},
\[
H(\s)= \sum_{\langle ij\rangle^*} U(\s_i,\s_j)\;,
\]
with the same pair interaction $U$ as in \rf{U}.  The variables $\s_i$
still take values in the set $E=\{0,1,\ldots,q\}$.  We then can
introduce the same definitions as in Section \ref{sec:WR}, only
replacing $\langle ij\rangle$ by $\langle ij\rangle^*$ at all proper
places. For simplicity, we refrain from adding $*$'s to the quantities 
so defined. 
In particular, we write again $\G(z,q)$ for the set of Gibbs measures
in the present square-shaped model. 

Given a configuration 
$\s\in\O$, we call a subset $S$ of $\Z^2$ an \emph{occupied 
$*$cluster} if $S$ is a maximal $*$connected subset of $\{i\in\Z^2: \s_i\ne 
0\}$. Such an $S$ may be visualized as the connected component in 
$\R^2$ of the union of all squares with centers in $S$.
We say that two disjoint square-shaped particles are 
\emph{contiguous} if they touch each other along at least one half of a side 
(i.e., if their centers have Euclidean distance $2$ or $\sqrt 5$), and 
are separated by two empty sites.
A set $S$ of pairwise disjoint square particles will be called
a \emph{contiguity sea} if it is maximal connected relative to the 
contiguity relation, and the union of their closures surrounds each 
bounded set.

Our result for this model is the following.
\begin{thm}\label{th:sqWR}
If $q$ exceeds some sufficiently large $q_0$, there exist a critical
activity $z_c(q)\in\;]\,q^{1/3}/3,3\,q^{1/3}\,[ $ and numbers 
$0<\e(q)<1/3$ with 
$\e(q)\to0$ as $q\ti$ such that the following hold:

{\rm (i) }For $z>z_c(q)$, there exist $q$ distinct translation 
invariant `colored' phases $\mu_a\in\G(z,q)$, $a\in\{1,\ldots,q\}$. 
Relative to $\mu_a$, there exists almost surely an occupied sea of 
color $a$ containing any given site with probability at least
$1-\e(q)$.

{\rm (ii) }For $q_0/q\le z<z_c(q)$, there exists a translation invariant
disordered phase $\mu_\dis\in\G(z,q)$ such that with probability $1$
all occupied $*$clusters are finite, independently and randomly colored, and 
surrounded by a contiguity sea. Moreover,
$\mu_\dis(\s_i\ne 0)<1/4+\e(q)$ for all $i$.

{\rm (iii) }At $z=z_c(q)$, a first-order phase transition occurs, in 
the sense that there exist $q+1$ distinct phases $\mu_\dis,\mu_1,\ldots, 
\mu_q\in\G(z_c(q),q)$ exhibiting the properties above.
\end{thm}

In short, we have the following phase diagram:
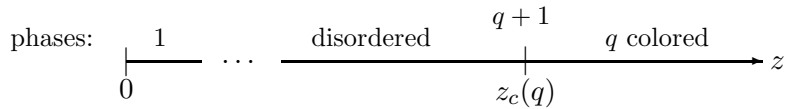
\begin{figure}[htb]
\begin{center}
\setlength{\unitlength}{1cm}
\begin{picture}(7.6,1)
\put(0.14,0){\line(1,0){1}}
\put(1.4,-0.024){$\ldots$}
\put(2.2,0){\vector(1,0){6.4}}
\put(8.7,-0.1){$z$}
\put(0,-0.31){$\mbox{$|$}\atop\mbox{0}$}
\put(5,-0.32){$\mbox{$|$}\atop\mbox{$z_c(q)$}$}
\put(-1.4,0.2){\small phases:}
\put(0.5,0.2){\small 1}
\put(2.6,0.2){\small disordered}
\put(6.5,0.2){\small $q$ colored}
\put(5,0.5){\small $q+1$}
\end{picture}
\end{center}
\end{figure}

\medskip\noindent
{\bf Remark \thethm } (1) The first-order transition at $z_c(q)$ manifests itself
thermodynamically by a jump of the particle density from a value 
close to $1/4$ to a value close to $1$.

\smallskip\noindent 
(2) The behavior of the square-shaped Widom--Rowlinson 
model differs from that of the standard, diamond-shaped Widom--Rowlinson 
model in that there are no staggered phases but instead only one
disordered phase showing not only color-disorder but also positional 
disorder. In fact, we expect that this disordered phase is the 
unique Gibbs measure for any $z<z_c(q)$, although this does not follow from 
our methods. (Just as in Remark \ref{th:WR} (5), one can 
see that for large $q$ the square-shaped Widom--Rowlinson model is 
related to the square-shaped hard-core lattice gas. It seems that 
the latter model does not exhibit a phase transition, but we are 
not aware of any proof.) 

\smallskip\noindent
(3) The first-order transition at $z_c(q)$ implies a percolation
transition from an empty sea in the disordered phase to an occupied
sea in the colored phases.  In spite of the supposed uniqueness of the
Gibbs measure in the whole range $[0,z_c(q)[$, this interval contains
a further percolation threshold, namely a critical value for the
existence of a contiguity sea.  Indeed, for sufficiently small $z$ one
can use disagreement percolation \cite{vdBM,GHM} to show the
uniqueness of the Gibbs measure and the existence of a sea of empty
plaquettes; the latter excludes the existence of a contiguity sea. 
However, it is far from obvious that a contiguity sea will set on
at a well-defined activity.  This is because neither the existence of
a contiguity sea is an increasing event, nor the site-random cluster
distributions (which can also be used in the present case) are
stochastically monotone in $z$, cf.\ Remark \ref{th:WR} (6).
We note that a similar percolation transition occurs also
in the square-shaped hard-core lattice gas \cite{FH}; this shows again
that the latter describes the limiting behavior of our model in this regime.

\subsection{A Widom--Rowlinson model with molecular hard core}
\label{sec:WRhc}

Another way of changing the geometry of the Widom--Rowlinson model is 
to introduce a molecular (i.e., color-independent) hard-core 
interaction between particles. In this section we will discuss a 
model variant of this kind.

As before, the underlying lattice is still the square lattice $\Z^2$, 
and the state space at each lattice site is the set $E=\{0,1,\ldots,q\}$.
The formal Hamiltonian is of the form
\be{Ham_hc}
H(\s)= \sum_{|i-j|=1} \Phi(\s_i,\s_j)+
\sum_{|i-j|=\sqrt 2} U(\s_i,\s_j)\;;
\ee
here $\Phi$, the nearest-neighbor molecular hard-core exclusion, is given 
by
\[
\Phi(\s_i,\s_j)=\left\{\ba{cl}\infty&\mbox{if }\s_i\s_j\ne 0 
\;,\\
0&\mbox{otherwise,}\ea
\right.
\]
and the next-nearest neighbor color repulsion $U$ is still defined by 
\rf{U}. 
The main effect of the molecular hard core is a richer high-density
phase diagram containing $2q$ phases with color-order \emph{and}
staggered positional order. The low-density regime is disordered both 
in the sense of color and position, as in the case of the square-shaped
Widom--Rowlinson model. The transition between these regimes is still of 
first order, though the positional order of the high-density phases is an 
impediment for this to occur. We do not repeat the definitions of
admissible configurations and of Gibbs measures, which are
straightforward.
\begin{thm}\label{th:WRhc}
If $q\ge q_0$ for a suitable $q_0$, there exist a 
threshold $z_c(q)\in\; ]\,q/18,18\,q\,[$ and numbers $0<\e(q)<1/5$ with 
$\e(q)\to 0$ as $q\ti$ such that the following hold:

{\rm (i) }For $z>z_c(q)$ there exist $2q$ distinct colored and 
staggered phases 
$\mu_{a,\even},\, \mu_{a,\odd}\in\G(z,q)$, 
$a\in\{1,\ldots,q\}$, which are invariant under even translations. 
Relative to $\mu_{a,\even}$ there exist almost surely both an even 
occupied $*$sea of color $a$ and an odd empty $*$sea, and any
two adjacent sites belong to these $*$seas with probability at least
$1-\e(q)$. $\mu_{a,\odd}$ is obtained from 
$\mu_{a,\even}$ by a one-step translation.

{\rm (ii) }For $q_0^{1/7}/q\le z<z_c(q)$, there exists a  
translation invariant
disordered phase $\mu_\dis\in\G(z,q)$ such that with probability $1$
all occupied $*$clusters are finite, independently colored with 
uniform distribution, and 
enclosed by a contiguity sea. Also,
$\mu_\dis(\s_i\ne 0)<1/4+\e(q)$ for all $i$.

{\rm (iii) }At $z=z_c(q)$, a first-order phase transition occurs, in 
the sense that there coexist $2q+1$ distinct phases 
$\mu_\dis,\,\mu_{1,\even},\,\mu_{1,\odd},\ldots,\, 
\mu_{q,\even},\, \mu_{q,\odd}\in\G(z_c(q),q)$,  with
the properties above.
\end{thm}

We summarize this theorem by the following phase diagram:

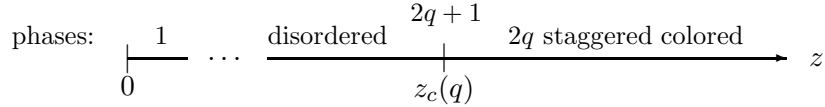
\begin{figure}[h!]
\begin{center}
\setlength{\unitlength}{1cm}
\begin{picture}(8.1,1)
\put(0.14,0){\line(1,0){0.8}}
\put(1.2,-0.024){$\ldots$}
\put(2.0,0){\vector(1,0){6.9}}
\put(9.2,-0.1){$z$}
\put(0,-0.31){$\mbox{$|$}\atop\mbox{0}$}
\put(3.9,-0.32){$\mbox{$|$}\atop\mbox{$z_c(q)$}$}
\put(-1.4,0.2){\small phases:}
\put(0.5,0.2){\small 1}
\put(2.0,0.2){\small disordered}
\put(5.2,0.2){\small $2q$ staggered colored}
\put(3.9,0.5){\small $2q+1$}
\end{picture}
\end{center}
\end{figure}

\medskip\noindent
{\bf Remark \thethm } (1) The first-order transition at $z_c(q)$ manifests itself
thermodynamically by a jump of the particle density from a value 
close to $1/4$ to a value close to $1/2$.

\smallskip\noindent 
(2) The Widom--Rowlinson model with molecular hard-core may be viewed
as a combination of a lattice gas of hard diamonds and the
square-shaped Widom--Rowlinson model.  Its high-density regime
inherits the staggered occupation pattern from the former, and the
color order from the latter.  The effect of colors is still strong
enough to produce a first-order transition, which is absent in the
pure hard-diamonds model \cite{LMNS}.  Just as in the case of Theorem
\ref{th:sqWR}, the low density regime is governed by the behavior of
the square-shaped hard-core lattice gas.  Indeed, it is not difficult
to develop a random-cluster representation of the model and to show
that its conditional probabilities converge to that of the
square-shaped lattice gas when $q\to\infty$ but $qz$ remains fixed;
cf.\ Remark \ref{th:WR} (5). Mutatis mutandis, the comments in Remark 
\ref{th:sqWR} apply here as well.

\smallskip\noindent 
(3) One may ask what happens if we interchange the r\^oles of $\Phi$ 
and $U$ in the Hamiltonian of Eq.\ \rf{Ham_hc}, i.e., if there is a 
molecular hard core between diagonal neighbors of distance $\sqrt 2$,
and a Widom--Rowlinson intercolor repulsion between nearest 
particles of distance $1$. In this case it is not hard to see that for 
any $q\ge 2$ and sufficiently large $z$ there exist four different 
phases with positional order but color disorder. One of these phases
(to be called the even vertical phase) has almost surely a sea of 
sites $i=(i_1,i_2)\in\Z^2$ which are occupied when $i_1$ is even, and 
empty when $i_1$ is odd. The other three phases are obtained
by translation and/or interchange of coordinates. However, it 
seems that in this case the geometry of interaction does not exhibit 
the properties leading to a 
first-order transition, so that the transition to the low density regime 
is of second order. We will return to 
this point at the end of Section \ref{sec:pf_hc}.

\subsection{A continuous-spin Widom--Rowlinson model}
\label{On-model}

The multitype Widom--Rowlinson lattice model may be viewed as a
diluted clock model for which each lattice site is either empty or
occupied by a particle with an orientation in the discrete group of
$q$'th roots of unity.  This suggests considering the following
plane-rotor model of oriented particles which may serve as a
simple model of a ferrofluid or liquid crystal.

Consider the state space $E=\{0\}\cup S^1$, equipped with the 
reference measure $\nu=\d_0+\l$, where $\l$ is normalized Haar
measure on the circle $S^1$. As before, the equality $\s_i=0$ means 
that site $i$ is empty, while $\s_i=a\in S^1$ says that $i$ is 
occupied by a particle with orientation $a$. The formal 
Hamiltonian is again given by \rf{Ham}, where the pair interaction $U$ 
is now defined by
\[
U(\s_i,\s_j)=\left\{\ba{cl}\infty&\mbox{if }\s_i, \s_j\in S^1,\ 
\s_i\cdot \s_j\ge \cos 2\pi\a \;,\\
0&\mbox{otherwise}\ea\right.
\]
for some angle $0<\a<1/4$.  This potential forces adjacent particles
to have nearly the same orientation.  The parameter $\a$ will play the
same r\^ole as $1/q$ did before.  The Gibbs distributions
$\mu_{\L,\eta}^{z,\a}$ in a finite region $\L\subset\Z^2$ with
boundary condition $\eta$ and activity $z>0$ are defined by their densities
with respect to the product measure $\nu^\L$, which are again given by
the right-hand side of equation \rf{Gibbs}.  We write $\G(z,\a)$ for
the associated set of Gibbs measures.  Since $U$ preserves the
$O(2)$-symmetry of particle orientations, the
Mermin--Wagner--Dobrushin--Shlosman theorem (cf.\ Theorem (9.20) of
\cite{Gii}) implies that each such Gibbs measure is invariant under
simultaneous rotations of particle orientations.
\begin{thm}\label{th:contspin}
If $\a$ is less than some sufficiently small $\a_0$,
there exist a critical activity $z_c(\a)\in\;]\,\a^{-2}/18,5\,\a^{-2}\,[$ 
and numbers $0<\e(\a)<1/3$ with $\e(\a)\to 0$ as $\a\to 0$ such that 
the following hold:

{\rm (i) }For $z>z_c(\a)$ there exists a dense `ordered' phase
$\mu_\ord\in\G(z,\a)$ exhibiting the translation invariance and
$O(2)$-symmetry of the model.  Relative to $\mu_\ord$, there exists
almost surely an occupied sea containing any fixed site with
probability at least $1-\e(\a)$ (and on which the orientations of
adjacent particles differ only by the angle $2\pi\a$).

{\rm (ii) }For $\a_0\inv\le z<z_c(\a)$, there exist two distinct
`staggered' phases $\mu_\even,\mu_\odd\in\G(z,\a)$ which are invariant
under particle rotations and even translations.  Almost surely with
respect to $\mu_\even$ there exist both an even occupied $*$sea and an
odd empty $*$sea, and any two adjacent sites belong to these $*$seas
with probability at least $1-\e(q)$.  In addition, all occupied
clusters are almost surely finite, and conditionally on their position
the distribution of orientations is invariant under simultaneous
rotations of all spins in a single occupied cluster.  $\mu_\odd$ is
obtained from $\mu_\even$ by a one-step translation.

{\rm (iii) }At $z=z_c(\a)$, a first-order phase transition occurs, in
the sense that there exist three distinct phases $\mu_\even, \mu_\odd,
\mu_\ord \in\G(z_c(\a),q)$ with the properties above.
\end{thm}
The theorem above shows that the present model behaves similarly to the
related finite-energy model considered in \cite{CKS2}.  Presumably
$\mu_\ord$ is the unique Gibbs measure for $z>z_c(\a)$, and there is a
second-order transition from the staggered regime to the low-activity
uniqueness regime.  We thus have the following phase diagram.
\begin{figure}[ht!]
\begin{center}
\setlength{\unitlength}{1cm}
\begin{picture}(7.6,1)
\put(0.14,0){\line(1,0){1}}
\put(1.4,-0.024){$\ldots$}
\put(2.2,0){\vector(1,0){6.4}}
\put(8.7,-0.1){$z$}
\put(0,-0.31){$\mbox{$|$}\atop\mbox{0}$}
\put(5,-0.32){$\mbox{$|$}\atop\mbox{$z_c(q)$}$}
\put(-1.4,0.2){\small phases:}
\put(0.5,0.2){\small 1}
\put(2.6,0.2){\small 2 staggered}
\put(6.5,0.2){\small ordered}
\put(5.3,0.5){$3$}
\end{picture}
\end{center}
\end{figure}

\medskip\noindent
{\bf Remark \thethm } \ The model above is a continuous-spin 
counterpart of the standard Widom--Rowlinson model considered in 
Section \ref{sec:WR}. It is rather straightforward to modify our 
techniques for investigating 
analogous continuous-spin variants of the square-shaped 
Widom--Rowlinson model and of the model with diagonal molecular hard 
core. In the first case, we obtain a phase diagram of the form

\begin{figure}[ht!]
\begin{center}
\setlength{\unitlength}{1cm}
\begin{picture}(7.6,1)
\put(0.14,0){\line(1,0){1}}
\put(1.4,-0.024){$\ldots$}
\put(2.2,0){\vector(1,0){6.4}}
\put(8.7,-0.1){$z$}
\put(0,-0.31){$\mbox{$|$}\atop\mbox{0}$}
\put(5,-0.32){$\mbox{$|$}\atop\mbox{$z_c(q)$}$}
\put(-1.4,0.2){\small phases:}
\put(0.5,0.2){\small 1}
\put(2.6,0.2){\small disordered}
\put(6.5,0.2){\small ordered}
\put(5.3,0.5){$2$}
\end{picture}
\end{center}
\end{figure}
\noindent
and in the second case we find

\begin{figure}[h!]
\begin{center}
\setlength{\unitlength}{1cm}
\begin{picture}(7.6,1)
\put(0.14,0){\line(1,0){1}}
\put(1.4,-0.024){$\ldots$}
\put(2.2,0){\vector(1,0){6.4}}
\put(8.7,-0.1){$z$}
\put(0,-0.31){$\mbox{$|$}\atop\mbox{0}$}
\put(4.1,-0.32){$\mbox{$|$}\atop\mbox{$z_c(q)$}$}
\put(-1.4,0.2){\small phases:}
\put(0.5,0.2){\small 1}
\put(2.4,0.2){\small disordered}
\put(5.2,0.2){\small  2 staggered ordered}
\put(4.5,0.5){$3$}
\end{picture}
\end{center}
\end{figure}
\noindent
The details are left to the reader.

\section{Proof of Theorem \ref{th:WR}}
\label{sec:pf_WR}

The proof of all four theorems follows the general scheme described in
Chapters 18 and 19 of Georgii \cite{Gii}, which is similar in spirit
to that of Dobrushin and Shlosman \cite{DS} and Koteck\'y and Shlosman
\cite{KS}.  This scheme consists of two parts: a model-specific
contour estimate implying percolation of ``good plaquettes'', and a
general part deducing from this percolation the first-order transition
and the properties of phases.  We describe the general part first and
defer the contour estimate to a second subsection.  Many of the
details presented here for the Widom--Rowlinson model carry over to
the other models, so that for the proofs of Theorems \ref{th:sqWR} to
\ref{th:contspin} we only need to indicate the necessary changes.  We
note that our arguments can easily be extended to the higher
dimensional lattices $\Z^d$ using either the ideas of Chapter 18 of
\cite{Gii} or those of \cite{Fuku}.

\subsection{Competition of staggered and ordered plaquettes}
\label{sec:compete}

We consider the standard plaquette $C=\{0,1\}^2$ in $\Z^2$ as well as
its translates $C+i$, $i\in\Z^2$.  Two plaquettes $C+i$ and $C+j$ will
be called adjacent if $|i-j|=1$, i.e., if $C+i$ and $C+j$ share a
side.  We are interested in plaquettes with a specified configuration
pattern.  Each such pattern will be specified by a subset $F$ of
$\O_C$, the set of admissible configurations in $C$.  For any such $F$
we define a random set $V(F)$ as follows.  Let $r_1$ and $r_2$ be the
reflections of $C$ in the vertical resp.\ horizontal line in the
middle of $C$, and $r^i= r_1^{i_1} r_2^{i_2}$ the reflection
associated to $i=(i_1,i_2) \in\Z^2$.  We then let
\be{V}
V(F): \s \to \{i\in\Z^2: r^i\s_{C+i}\in F\}
\ee
be the mapping associating with each $\s\in\O$ the set of
plaquettes on which $\s$ shows the pattern specified by $F$.  (The
reflections $r^i$ need to be introduced for reasons of consistency:
they guarantee that two adjacent plaquettes may both belong to $V(F)$
even when $F$ is not reflection invariant, as e.g.\ the sets $G_\even
$ and $G_\odd$ below.)

We are interested in the case when $F$ is one of the following sets
of `good' configurations on $C$. These sets are distinguished
according to their occupation pattern. Describing a configuration on 
$C$ by a $2\times 2$ matrix in the obvious way, we define 
\bit 
\item $G_\st = G_\even \cup G_\odd \equiv
\{{0\;b \choose a\;0}: 1\le a,b\le q\}\cup 
\{{a\;0 \choose 0\;b}: 1\le a,b\le q\}$, the set of all
{\em staggered\/} configurations with `diagonal occupations'.

\item $G_\ord=\bigcup_{1\le a\le q} G_a\equiv 
\bigcup_{1\le a\le q}\{{a\;a \choose a\;a}\}$, the set of all
fully {\em ordered\/} configurations with four particles 
of the same color.

\item $G=G_\st\cup G_\ord$, the set of all good configurations.
\eit
Our first objective is to establish percolation of good plaquettes,
i.e., of plaquettes in which the configuration is good; the other
plaquettes will be called bad.  We want to establish this kind of
percolation for suitable Gibbs measures \emph{uniformly in the
activity $z$} (provided $z$ is not too small).  A suitable class of
Gibbs measures is that obtained by infinite-volume limits with
periodic boundary conditions.

For any integer $L\ge1$ we consider the rectangular box 
\be{LaL}
\L_L=\{-12\,L+1,\ldots,12\,L\}\times \{-14\,L+1,\ldots,14\,L\}
\ee
in $\Z^2$ of size $v(L)=24\,L\times 28\, L$.  (The reason for this
particular choice will become clear in the proofs of Lemmas \ref{B3}
and \ref{B2}.)  We write $\mu_{L,\per}^{z,q}$ for the Gibbs
distribution in $\L_L$ with parameters $z,q$ and periodic boundary
condition, and $\G_\per(z,q)$ for the set of all limiting measures of
$\mu_{L,\per}^{z,q}$ as $L\ti$ (relative to the weak topology of
measures).  The basic result is the following {\em contour estimate\/}
which shows that bad plaquettes have only a small chance to occur.
%
\begin{prop}\label{contour-est} For any $\d>0$ there exists a 
number $q_0\in\N$ such that 
\be{cont-est}
\mu(\D \cap V(G)=\emptyset) \le \d^{|\D |}
\ee
whenever $q\ge q_0$, $zq\ge q_0$, $\mu\in \G_\per(z,q)$, and $\D \subset 
\Z^2$ is finite.
\end{prop}
In the above, $\{\D \cap V(G)=\emptyset\}$ is a short-hand for the event
consisting of all $\s$ for which all plaquettes $C+i$, $i\in \D $, are bad;
similar abbreviations will also be used below.

The proof of the proposition takes advantage of reflection positivity and 
the chessboard estimate, cf.\ Corollary (17.17) of \cite{Gii}, and is 
deferred to the next section. We mention here only that $q_0$ is 
chosen so large that
\be{q_0}
\d(q)\equiv q^{-1/56}+q^{-1/12}+q^{-1/4}+q^{-1/2}\le\d
\ee
when $q\ge q_0$. It will be essential in the following that the 
contour estimate is uniform for $z\ge q_0/q$.  

As an immediate consequence of the contour estimate we obtain the
existence of a \emph{sea of good plaquettes}. We will say that a set of 
plaquettes forms a sea if the set of their left lower corners is 
connected and surrounds each finite set. It is then evident that the 
existence of a sea of completely occupied plaquettes implies the 
existence of an occupied sea; likewise, the existence of a sea of 
plaquettes which are occupied on their even points implies the 
existence of an even occupied $*$sea. In this way, the concept of a 
sea of plaquettes is general enough to include all concepts of seas 
introduced in Section \ref{results}. 
Specifically, for any $F\subset \O_C$ we
define $\sea(F)$ as the largest sea in $V(F)$ whenever $V(F)$ contains
a sea, and let $\sea(F)=\emptyset$ otherwise. 

For $z\ge q_0/q$ we
write $\bar\G_\per(z,q)$ for the set of all accumulation points (in
the weak topology) of measures $\mu_n\in\G_\per(z_n,q)$ with $z_n\to
z$, $z_n\ge q_0/q$. The graph of the correspondence $z\to 
\bar\G_\per(z,q)$ is closed; this will be needed in the proof of 
property (A2) below.
%
\begin{prop}\label{sea} For any $\e>0$ there exists a 
number $q_0\in\N$ such that 
\[
\mu(0\in \sea(G))\ge 1-\e
\]
whenever $\mu\in\bar\G_\per(z,q)$, $q\ge q_0$ and $z\ge q_0/q$.
\end{prop}
\proof Note first that the contour estimate \rf{cont-est} involves
only local events and therefore extends immediately to all $\mu\in
\bar\G_\per(z,q)$.  The statement then follows directly from
Proposition \ref{contour-est} together with Lemmas (18.14) and (18.16)
of \cite{Gii}.  The number $\d$ has to be chosen so small that
$4\d(1-5\d)^{-2}\le\e\;$.  $\Box$

\medskip\noindent
What is the advantage of having a sea of good plaquettes?
The key property is that the sets $G_\st$ and $G_\ord$ have disjoint 
side-projections. That is, writing $b=\{(0,0),(1,0)\}$ for the two points
on the bottom side of $C$ we have
\[
\s\in G_\st, \s'\in G_\ord \Rightarrow \s_b\ne\s_b'\;,
\]
and similarly for the other sides of $C$.  As a consequence, if two
adjacent plaquettes are good then they are both of the same type,
either staggered or ordered.  Therefore each sea of good plaquettes is
either a sea of staggered plaquettes, or a sea of ordered plaquettes. 
Hence
\[
\{\sea(G)\ne\emptyset\}=\{\sea(G_\st)\ne\emptyset\}\cup 
\{\sea(G_\ord)\ne\emptyset\}\;,
\]
and the two sets on the right-hand side are disjoint.
Moreover, the sets $G_\even$ and $G_\odd$ also have disjoint 
side-projections, and so do the sets $G_a$, $1\le a\le q$. 
Therefore, the event $\{\sea(G_\st)\ne\emptyset\}$ splits into the 
two disjoint subevents $\{\sea(G_\even)\ne\emptyset\}$ and 
$\{\sea(G_\odd)\ne\emptyset\}$, and $\{\sea(G_\ord)\ne\emptyset\}$ 
splits off into the disjoint subevents $\{\sea(G_a)\ne\emptyset\}$,
$1\le a\le q$. In other words, each sea of good plaquettes has a 
characteristic occupation pattern or color corresponding to a
particular phase, and we only need to identify the activity regimes
for which the different phases do occur.
 
To this end we fix any $\e>0$.  We will need later that $\e<1/6$.  As
in the proof of Proposition \ref{sea}, we choose some $0<\d<1/25$ such
that $4\d(1-5\d)^{-2}\le\e\;$, and we let $q_0$ be so large that
condition \rf{q_0} holds for all $q\ge q_0$.  For such $q_0$ and $q$
we consider the two activity domains
\[
A_\st = \Big\{z\ge q_0/q:  \mu(0\in V(G_\st))\ge \mu(0\in V(G_\ord))
\mbox{ for some }\mu\in \bar\G_\per(z,q)\Big\} 
\]
and
\[
A_\ord = \Big\{z\ge q_0/q:  \mu(0\in V(G_\ord))\ge \mu(0\in V(G_\st))
\mbox{ for some }\mu\in \bar\G_\per(z,q)\Big\}\;.
\]
Our next result shows that these sets describe the regimes in which
staggered resp.\ ordered phases exist.  The \emph{mean particle
density} $\r(\mu)$ of a measure $\mu$ which is periodic under
translations with period 2 is defined by
\be{rho}
\r(\mu)= \mu(N_C)/|C|\;;
\ee
recall that $N_C$ is the number of particles in $C$.
%
\begin{prop}\label{phases}   
{\rm(a) } For each $z\in A_\st$ there exist two `staggered' Gibbs
measures $\mu_\even,\mu_\odd\in\G(z,q)$ invariant under even
translations of $\Z^2$ and permutations of particle colors. 
$\mu_\even$-almost surely we have $\sea(G_\even)\ne\emptyset$, and all
occupied clusters are finite and have independently distributed random
colors.  In addition, $\mu_\even(0\in\sea(G_\even))\ge1-2\e$, and in
particular $\r(\mu_\even)\le \frac12+\e$.  $\mu_\odd$ has the
analogous properties.

{\rm(b) } For each $z\in A_\ord$ there exist $q$ `colored' translation 
invariant Gibbs measures $\mu_a\in\G(z,q)$, $a\in\{1,\ldots,q\}$ . 
Each $\mu_a$ satisfies $\mu_a(\sea(G_a)\ne\emptyset)=1$,
$\mu_a(0\in\sea(G_a))\ge1-2\e$, and in particular
has mean particle density $\r(\mu_a)\ge 1-2\e$.
\end{prop}
\proof (a) Let $z\in A_\st$ be given and $\mu\in\bar\G_\per(z,q)$ 
be such that $\mu(0\in V(G_\st))\ge \mu(0\in V(G_\ord))$. Then
$\mu(0\in V(G_\ord))\le 1/2$ and therefore
\bea
\mu\Big(0\in\sea(G_\st)\Big)&=&
\mu\Big(0\in\sea(G),\, 0 \not\in V(G_\ord)\Big)\\
&\ge&1-\e -\frac12= \frac12 -\e>0\;.
\eea
But $G_\st$ splits into 
the two parts $G_\even$ and $G_\odd$ which are related to each other 
by the reflection in the line $\{x_1=1/2\}$, and $\mu$ is invariant 
under this reflection. Hence 
\[
p\equiv \mu\Big(0\in\sea(G_\even)\Big)= \mu\Big(0\in\sea(G_\odd)\Big)\ge 
\frac12\left(\frac12-\e\right)>0\;.
\] 
We can therefore define the conditional probabilities
$\mu_\even=\mu(\dot|\sea(G_\even)\ne\emptyset)$ and
$\mu_\odd=\mu(\dot|\sea(G_\odd)\ne\emptyset)$. Since the events in 
the conditions are tail measurable, these measures belong to $\G(z,q)$.
It is clear that these conditional probabilities inherit all common 
invariance properties of $\mu$ and the conditioning events.
Moreover, we find
\bea
\mu\Big(\sea(G_\even)\ne\emptyset\Big)&=&\frac12\, 
\mu\Big(\sea(G_\st)\ne\emptyset\Big)\\
&\le& 
\frac12\,\mu\Big(0\in\sea(G_\st)\Big)+\frac12\,\mu\Big(0\not\in\sea(G)\Big)\\
&\le& p +\frac{\e}{2}\;,
\eea
and therefore 
$$\mu_\even\Big(0\in\sea(G_\even)\Big)\ge\frac{p}{p+\e/2}\ge 1-2\e\;.$$
In particular, it follows that 
\[
\r(\mu_\even) \le \frac12\; \mu_\even\Big(0\in V(G_\even)\Big) + 
\mu_\even\Big(0\not\in V(G_\even)\Big)\le \frac12 + \e\;.
\]
Finally, we show that $\mu_\even$-almost surely all occupied clusters
are finite, and their colors are conditionally independent and
uniformly distributed when all particle positions are fixed.  Indeed,
since $\mu_\even(\sea(G_\even)\ne\emptyset)=1$ there exists
$\mu_\even$-almost surely an odd empty $*$sea.  This means that any
box $\D$ is almost surely surrounded by an empty $*$circuit.  On the
one hand, this shows that all occupied clusters must be finite almost
surely.  On the other hand, for any $\eta>0$ we can find a box
$\D'\supset\D$ containing an empty $*$circuit around $\D$ with
probability at least $1-\eta$.  Let $\Gamma$ be the largest set with
$\D\subset\Gamma\subset\D'$ such that there are no particles on its
outer boundary $\partial\Gamma$; if no such set exists we set
$\Gamma=\emptyset$.  The events $\{\Gamma=\L\}$ then depend only on
the configuration in $\Z^2\setminus \L$. By the strong Markov
property of $\mu_\even$, we conclude that on 
$\{\Gamma\ne\emptyset\}$ the distribution of colors of the occupied
clusters meeting $\D$ is governed by the Gibbs distribution in
$\Gamma$ with empty boundary condition.  The symmetry properties of
the latter thus imply that these colors are conditionally independent
and uniformly distributed.  Letting $\eta\to 0$ and $\D\uparrow\Z^2$ we
find that this statement holds in fact for all occupied clusters.

By construction, $\mu_\odd$ is obtained from $\mu_\even$ by a 
one-step translation, and thus has the analogous properties. 

(b) The proof of this part is quite similar. Pick any 
$z\in A_\ord$ and $\mu\in\bar\G_\per(z,q)$ such that 
$\mu(0\in V(G_\ord)|0\in V(G))\ge 1/2$. Since $\mu$ is invariant under 
permutations of colors it then follows in the same way that
\[
p\equiv \mu\Big(0\in\sea(G_a)\Big)\ge 
\frac1q\left(\frac12-\e\right)>0\;,
\]
so that we can define the conditional probabilities $\mu_a=
\mu(\dot|\sea(G_a)\ne\emptyset)\in\G(z,q)$, $a\in\{1,\ldots,q\}$.
Also, 
\[
\mu\Big(\sea(G_a)\ne\emptyset\Big)
\le 
\frac1q\,\mu\Big(0\in\sea(G_\ord)\Big)+\frac1q\,\mu\Big(0\not\in\sea(G)\Big)
\le p +\frac{\e}{q}\;,
\]
whence $\mu_a(0\in\sea(G_a))\ge{p}/({p+\e/q})\ge 1-2\e$ and
$
\r(\mu_a)\ge \mu_a(0\in V(G_a))\ge 1-2\e
$. 
$\Box$

\medskip
According to the preceding proposition, Theorem \ref{th:WR} will be
proved once we have shown that there exists a critical activity
$z_c(q)\in \;]\,q/5,5q\,[$ such that $A_\st =[q_0/q,z_c(q)]$ and
$A_\ord=[z_c(q),\infty[$.  To this end we will establish the following
items:
\newcounter{num}
\begin{list}{(A\arabic{num}) }{\usecounter{num}
\setlength{\leftmargin}{9ex} \setlength{\labelsep}{2ex}}
\item[(A1)] $A_\st\cup A_\ord= [q_0/q,\infty[\;$.
\item[(A2)] $A_\st$ and $A_\ord$ are closed.
\item[(A3)] $A_\ord\cap [q_0/q,{q/5}] \;=\emptyset\;$.
\item[(A4)] $A_\st\cap [5 q,\infty[ \;=\emptyset\;$.
\item[(A5)] $|A_\st\cap A_\ord|\le 1\;$.
\end{list}
Statement (A1) follows trivially from the definitions of $A_\st$ and $A_\ord$.
Assertion (A2) is also obvious because these definitions involve only 
local events, and the graph of the correspondence
$z\to\bar\G_\per(z,q)$ is closed by definition. 

Property (A3) corresponds to the discovery of Runnels and Lebowitz 
\cite{RL} that staggered phases do exist in a nontrivial activity regime,
and follows directly from the next result.
%
\begin{lem}\label{stag} 
For $z\le q/5$ and $\mu\in \bar\G_\per(z,q)$ we have
\[
\mu\Big(0\in V(G_\ord)\Big|0\in V(G)\Big) <1/2\;.
\]
\end{lem}
\proof Consider the Gibbs distribution $\mu_{L,\per}^{z,q}$ in 
the box $\L_L$ with periodic boundary condition, and let 
\be{Gord}
G_{\ord,L}=\Big\{\s\in \O_{L,\per}: \s_{C(i)}\in G_\ord 
\mbox{ for all }i\in\L_L\Big\}\,;
\ee
here we write $\O_{L,\per}$ for the set of admissible configurations 
in the torus $\L_L$ (including nearest-neighbor bonds between the left 
and the right sides as well as between the top and bottom 
sides of $\L_L$),
and $C(i)$ for the image $C+i \mbox{ mod }\L_L$ of $C$ under the 
periodic shift of $\L_L$ by $i$.
(As $G_\ord$ is reflection-symmetric, we can omit the 
reflections $r^i$ which appear in \rf{V}.)
The chessboard estimate (cf.\ Corollary (17.17) of \cite{Gii})
then implies that 
\[
\mu_{L,\per}^{z,q}\Big(0\in V(G_\ord)\Big) \le 
\mu_{L,\per}^{z,q}(G_{\ord,L})^{1/v(L)}\;.
\]
We compare the latter probability with that of the event
\be{Geven}
G_{\even,L}=\Big\{\s\in \O_{L,\per}: r^i\s_{C(i)}\in G_\even \mbox{ for 
all }i\in\L_L\Big\}\;.
\ee 
This gives
\[
\mu_{L,\per}^{z,q}(G_{\ord,L})\le 
\mu_{L,\per}^{z,q}(G_{\ord,L})\Big/\mu_{L,\per}^{z,q}(G_{\even,L})
= z^{v(L)}q \;z^{-v(L)/2}q^{-v(L)/2}
\]
because $G_{\ord,L}$ contains only the $q$ distinct close packed 
monochromatic configurations, while for $\s\in G_{\even,L}$ the 
$v(L)/2$ particles can have independent colors. Taking the $v(L)$'th 
root and letting $L\ti$ we find for $\mu\in\bar\G_\per(z,q)$
\[
\mu\Big(0\in V(G_\ord)\Big)\le (z/q)^{1/2}\le 5^{-1/2}<  
(1-\d)/2\;.
\]
The last inequality comes from the choice of $\d$.
Since $\mu(0\in V(G))\ge 1-\d$ by Proposition \ref{contour-est}, 
the lemma follows.
$\Box$

\medskip
Assertion (A4) corresponds to the well-known fact that $q$ ordered 
phases exist when the activity is large. For $q=2$ this was already 
shown by Lebowitz and Gallavotti \cite{LG}, and for arbitrary $q$ 
by Runnels and Lebowitz 
\cite{RL}. This is again a simple application of the chessboard 
estimate.
%
\begin{lem}\label{order} 
For $z\ge 5q$ and $\mu\in \bar\G_\per(z,q)$ we have
\[
\mu\Big(0\in V(G_\st)\Big|0\in V(G)\Big) <1/2\;.
\]
\end{lem}
\proof Let $G_{\ord,L}$ be as in \rf{Gord}, and define $G_{\st,L}$ 
analogously. By the chessboard estimate we find
\bea
\mu_{L,\per}^{z,q}\Big(0\in V(G_\st)\Big) &\le& 
\mu_{L,\per}^{z,q}(G_{\st,L})^{1/v(L)}\\
&\le& \Big(
\mu_{L,\per}^{z,q}(G_{\st,L})\Big/\mu_{L,\per}^{z,q}(G_{\ord,L})\Big)^{1/v(L)}\\
&\le& 2^{1/v(L)}z^{1/2}q^{1/2} \;z^{-1}q^{-1/v(L)}
\eea
because $G_{\st,L}=G_{\even,L}\cup G_{\odd,L}$ contains 
$2\,q^{v(L)/2}$ distinct 
configurations of particle density $1/2$. We can now complete the 
argument as in the preceding proof.
$\Box$

\medskip
For the proof of (A5) we will use a thermodynamic argument, namely the 
convexity of the pressure as a function of $\log z$. 
For any translation invariant probability measure $\mu$ on $\O$ we
consider the {\em entropy per volume} 
\[
s(\mu)=\lim_{|\L|\to\infty} |\L|^{-1}\, S(\mu_\L)\;.
\]
Here we write $\mu_\L$ for the restriction of $\mu$ to 
$\O_\L$, 
\[
S(\mu_\L)=-\sum_{\s\in\O_\L}\mu_\L(\s) \,\log \mu_\L(\s)
\]
is the entropy of $\mu_\L$, and the notation $|\L|\to\infty$ means that 
$\L$ runs through a 
specified increasing sequence of square boxes;
for the existence of $s(\mu)$ we refer to \cite{Gii,Rue}.

We  define the thermodynamic \emph{pressure} by
\be{p}
P(\log z)=\max_\mu \Big[\r(\mu)\,\log z +s(\mu)\Big]\;;
\ee
the maximum extends over all 
translation invariant probability measures $\mu$ on $\O$,
and $\r(\mu)=\mu(\s_0\ne 0)$ is the associated mean particle density,
cf. \rf{rho}.
(Since $\O$ is defined as the set of all admissible 
configurations, the hard-core intercolor repulsion is taken into 
account automatically.)
By definition, $P$ is a convex function of $\log z$, and the 
variational principle (see Theorems 4.2 and 3.12 of \cite{Rue}) 
asserts that the maximum in \rf{p} is attained precisely on $\G_\Th(z,q)$,
the set of all translation invariaion invariant elements of $\G(z,q)$.
By standard arguments (cf.\ Remark (16.6) and Corollary (16.15) of 
\cite{Gii}) it follows that $P$ is strictly convex, and 
\be{tangent}
P_-'(\log z) \le \r(\mu) \le P_+'(\log z)\mbox{ for all 
}\mu\in\G_\Th(z,q)\;;
\ee
here we write $P_-'$ and $P_+'$ for the left-hand resp.\ 
right-hand derivative of $P$.
By strict convexity, $P_-'$ and $P_+'$ are strictly
increasing and almost everywhere identical. Assertion (A5) thus 
follows from the lemma below.
%
\begin{lem}\label{z_c} 
For each $z\in A_\st\cap A_\ord$ we have
$P_-'(\log z) \le 2/3 \le P_+'(\log z)$.
\end{lem}
\proof This has already been shown essentially in 
Proposition \ref{phases}. Pick any $z\in A_\st\cap A_\ord$, and 
let $\mu\in\bar\G_\per(z,q)$ be as in the 
proof of Proposition \ref{phases}(a). 
Consider the conditional probability 
$\mu_\st=\mu(\dot|\sea(G_\st)\ne\emptyset)=\frac12\,\mu_\even+
\frac12\,\mu_\odd$. By the arguments there, $\mu_\st$ is 
well-defined, belongs to $\G_\Th(z,q)$, and satisfies 
$\r(\mu_\st)\le 1/2 +\e < 2/3$. On the other hand, the measures
$\mu_a$ constructed in Proposition \ref{phases}(b) also belong to
$\G_\Th(z,q)$ and satisfy $\r(\mu_a)\ge 1-2\e>2/3$. The lemma thus 
follows from \rf{tangent}.
$\Box$

\medskip
We can now complete the proof of Theorem \ref{th:WR}.
Properties (A1) to (A4) together imply that $A_\st\cap A_\ord 
\ne\emptyset$. This is because the interval $[q_0/q,\infty[$ is 
connected and therefore cannot be the union of two disjoint non-empty 
closed sets. Combining this with (A5) we find
that $A_\st\cap A_\ord$ consists of a unique value 
$z_c(q)$. In particular, $A_\ord$ cannot contain any value 
$z<z_c(q)$ because the infimum of such $z$'s would belong to $A_\st\cap 
A_\ord$; likewise, $A_\st$  does not contain any value 
$z>z_c(q)$. Hence $A_\st =[q_0/q,z_c(q)]$ and 
$A_\ord=[z_c(q),\infty[$, and Theorem \ref{th:WR} follows from 
Proposition \ref{phases}.

\subsection{Contour estimates}
\label{sec:cont_est}

In this subsection we will prove Proposition \ref{contour-est}.
Consider the set $\O_C$ of all
admissible configurations in $C$, and the set $B=\O_C\setminus G$ of all 
bad configurations in $C$. We split $B$ into the following subsets 
which are distinguished by their occupation pattern:
\bit 
\item $B_0 =\{{0\;0 \choose 0\;0}\}$, 
the singleton consisting of the empty configuration in $C$.

\item $B_1 =
\{{0\;0 \choose a\;0}, {0\;0 \choose 0\;a}, {0\;a \choose 0\;0},
{a\;0 \choose 0\;0}: 1\le a\le q\}$, 
the set of all configurations with a single particle in $C$.

\item $B_{2}=\{{a\;a \choose 0\;0},{a\;0 \choose a\;0}, {0\;0 \choose
a\;a},{0\;a \choose 0\;a} : 1\le a\le q\}$, the set of admissible
configurations for which one side of $C$ is occupied, and the other
side is empty.

\item $B_3 
=\{{a\;a \choose 0\;a}, {a\;a \choose a\;0}\}, {a\;0 \choose a\;a},
{0\;a \choose a\;a}: 1\le a\le q\}$, the set of all
admissible configurations with three particles in $C$.
\eit
We then clearly have
$B = \bigcup_{k=0}^3 B_k$.
The four different kinds of ``badness'' of a plaquette
will be treated separately in the three lemmas below. We start with the most 
interesting case of plaquettes with three particles.

For any $L\ge1$ and $k\in \{0,\ldots,3\}$ let 
\[
B_{k,L}=\{\s\in \O_{L,\per}: \s_{C(i)}\in  B_k  \mbox{ for 
all }i\in\L_L\}\;,
\]
where $C(i)$ is as in \rf{Gord}. Consider the quantities 
$p_{k,L}^{z,q}= \mu_{L,\per}^{z,q}(B_{k,L})^{1/v(L)}$ and 
$p_k^{z,q}=\limsup_{L\ti}p_{k,L}^{z,q}$.
%
\begin{lem}\label{B3} 
$p_3^{z,q}\le q^{-1/56}$ for all $z>0$ and $q\in\N$.
\end{lem}
\proof 
Fix any integer $L\ge 1$ and consider the set $B_{3,L}$ of 
configurations $\s$ in $\L_L$ having a single empty site in each 
plaquette. We claim that $|B_{3,L}|< q\;2^{14L+2}$. First of all, for 
each $\s\in B_{3,L}$
the occupied sites in $\L_L$ form a connected set, so that all 
particles have the same color. Thus there are only $q$ possible colorings, 
and we only need to count the possible occupation patterns for $\s\in B_{3,L}$.
It is easy to see that the plaquettes $C(i)$ with $\s(i)=0$ form a 
partition of $\L_L$. For each such partition, the 
plaquettes are either arranged in rows or in 
columns. In the first case, each row is determined by its parity (even or odd), 
namely the parity of $i_1$ for each $C(i)$ in this row; 
likewise, in the second case each column is determined by its parity. 
We can therefore count all such partitions as follows. There are 4 
possibilities of choosing the plaquette containing the origin. If 
this plaquette is fixed, there are no more than 
$2^{14L-1}$ possibilities of arranging all plaquettes in rows and choosing 
the parity of each row. Similarly, there are at
most $2^{12L-1}$ possibilities of arranging the plaquettes in 
columns. The number of such partitions is therefore no larger than
$4(2^{14L-1}+2^{12L-1})$, and the claim follows.

To estimate $\mu_{L,\per}^{z,q}(B_{3,L})$ we will rearrange the positions of 
all particles so that many different colors become possible. More 
precisely, we divide $\L_L$ into $(3L)(4L)$ rectangular cells $\D(j)$ of 
size $8\times7$. Let $\D_0(j)$ be the 
rectangular cell of size $7\times6$ situated in the left lower corner 
of $\D(j)$, and consider the set
\[
F_{3,L}=\{\s\in\O_{L,\per}: \s\ne 0\mbox{ on }\D_0(j), \;
\s\equiv 0\mbox{ on }\D(j)\setminus\D_0(j)\mbox{ for all }j\}\;.
\]
Since $|\D(j)\setminus\D_0(j)|=8\cdot 7- 7\cdot 6= |\D(j)|/4$ for 
all $j$, each $\s\in F_{3,L}$ has particle number $3\,v(L)/4$, just as 
the configurations in $B_{3,L}$. As the colors of the particles in 
the blocks $\D_0(j)$ can be chosen independently, we have
$|F_{3,L}|=q^{12L^2}=q^{v(L)/56}$. (The above construction, together 
with a similar construction in the proof of the next lemma, explains 
our choice of the rectangle $\L_L$.)
Now we can write
\[
\mu_{L,\per}^{z,q}(B_{3,L}) \le 
\frac{\mu_{L,\per}^{z,q}(B_{3,L})}{\mu_{L,\per}^{z,q}(F_{3,L})}
= \frac{|B_{3,L}|}{|F_{3,L}|} \le 2^{14L+2}\; q^{1-v(L)/56}\;.
\]
The proof is completed by taking the $v(L)$'th root and letting $L\ti$.
$\Box$

\medskip
Next we estimate the probability of plaquettes with two adjacent 
particles at one side of $C$.
%
\begin{lem}\label{B2} 
$p_2^{z,q}\le q^{-1/12}$ for all $z>0$ and $q\in\N$.
\end{lem}
\proof Fix any $L\ge1$, and let $\s\in B_{2,L}$. Then the particles 
are either arranged in alternating occupied and empty rows, or in 
alternating occupied and empty columns. The colors in all rows resp.\ 
columns can be chosen independently of each other. Hence
$|B_{2,L}|= 2(q^{14L}+q^{12L})\le 4\,q^{14L}$. Moreover, each $\s\in 
B_{2,L}$ has particle number $v(L)/2$. As in the last proof, we 
construct a set $F_{2,L}$ of configurations with the same particle number but 
larger color entropy as follows.

We partition $\L_L$ into $(8L)(7L)$ rectangular cells $\D(j)$ of 
size $3\times4$, and let $\D_0(j)$ be the 
rectangular cell of size $2\times3$ in the left lower corner 
of $\D(j)$. We then define
\[
F_{2,L}=\{\s\in\O_{L,\per}: \s\ne 0\mbox{ on }\D_0(j), \;
\s\equiv 0\mbox{ on }\D(j)\setminus\D_0(j)\mbox{ for all }j\}\;.
\]
Since $|\D(j)\setminus\D_0(j)|=3\cdot 4- 2\cdot 3= |\D(j)|/2$ for 
all $j$, each $\s\in F_{2,L}$ has particle number $v(L)/2$. 
As the particle colors in 
the blocks $\D_0(j)$ can be chosen independently, we have
$|F_{3,L}|=q^{56\,L^2}=q^{v(L)/12}$. As in the last proof, we thus 
find
\[
\mu_{L,\per}^{z,q}(B_{2,L}) \le 
\frac{\mu_{L,\per}^{z,q}(B_{2,L})}{\mu_{L,\per}^{z,q}(F_{2,L})}
= \frac{|B_{2,L}|}{|F_{2,L}|} \le 4\; q^{14L-v(L)/12}\;.
\]
Taking the $v(L)$'th root and letting $L\ti$ we obtain the result. 
$\Box$

\medskip
Finally we consider the probability of `diluted' plaquettes with a single 
or no particle. 
%
\begin{lem}\label{B10} 
$p_0^{z,q}\le (zq)^{-1/2}$ and $p_1^{z,q}\le (zq)^{-1/4}$ for all 
$z>0$, $q\in\N$.
\end{lem}
\proof We consider first the case of no particle.
For each $L\ge1$ we can write
\[
\mu_{L,\per}^{z,q}(B_{0,L}) \le 
\frac{\mu_{L,\per}^{z,q}(B_{0,L})}{\mu_{L,\per}^{z,q}(G_{\even,L})}
= \frac{1}{z^{v(L)/2}q^{v(L)/2}}\;,
\]
where $G_{\even,L}$ is defined by \rf{Geven}. 
The identity follows from the facts that 
$B_{0,L}$ contains only the empty configuration, whereas
each configuration in $G_{\even,L}$ consists of $v(L)/2$ 
particles with arbitrary colors. The first result is thus obvious.

Turning to the case of a single particle per plaquette, we note that 
each $\s\in B_{1,L}$ consists of $v(L)/4$ 
particles with arbitrary colors, and
there are no more than $2^{14L+2}$ distinct occupation patterns for 
these particles; the latter follows as in the proof of Lemma \ref{B3}
(by interchanging empty and occupied sites). Hence
\[
\mu_{L,\per}^{z,q}(B_{1,L}) \le 
\frac{\mu_{L,\per}^{z,q}(B_{1,L})}{\mu_{L,\per}^{z,q}(G_{\even,L})}
\le \frac{2^{14L+2} z^{v(L)/4}q^{v(L)/4}}{z^{v(L)/2}q^{v(L)/2}}\;,
\]
and the second result follows by taking the $v(L)$'th root and letting 
$L\ti$.
$\Box$

\medskip\noindent
{\em Proof of Proposition \ref{contour-est}. }Let $\mu\in 
\G_\per(z,q)$ and a finite $\D \subset\Z^2$ be given. Then we can write
\bea
\mu(\D \cap V(G)=\emptyset) &=& \sum_{\g:\D \to\{0,\ldots,3\}} 
\mu(\s: \s_{C+i}\in
B_{\g(i)} \mbox{ for all }i\in \D )\\
&\le& \sum_{\g:\D \to\{0,\ldots,3\}} \limsup_{L\ti} 
\mu_{L,\per}^{z,q}(\s: \s_{C(i)}\in
B_{\g(i)} \mbox{ for all }i\in \D )\\
&\le&\sum_{\g:\D \to\{0,\ldots,3\}} \limsup_{L\ti} \prod_{i\in \D } 
p_{\g(i),L}^{z,q}\\
&\le& \bigg( \sum_{k=0}^3 p_{k}^{z,q} \bigg)^{|\D |}
\;.
\eea
In the third step we have used the chessboard estimate, see Corollary (17.17) 
of \cite{Gii}. Inserting the estimates of Lemmas \ref{B3}, \ref{B2} 
and \ref{B10} and choosing $q_0$ as in \rf{q_0} we get the result.
$\Box$

\section{Proof of Theorem \ref{th:sqWR}}
\label{pf:sqWR}

Here we indicate how the proof of Theorem \ref{th:WR} can be adapted 
to obtain Theorem \ref{th:sqWR}. First of all, the different geometry 
of the present model leads to a new 
classification of good and bad plaquettes: the ordered configurations 
in $G_\ord$ are still good, but the (former good) configurations in 
$G_\st$ are now bad and will be denoted by $B_\st$, while the 
configurations in $B_1$ are now good, and we set $G_\dis=B_1$.

We first need an analog of the contour estimate, 
Proposition \ref{contour-est}. Remarkably, the estimates of Lemmas 
\ref{B3} and \ref{B2} carry over without any change. To deal with $B_\st$
we can proceed exactly as in Lemma \ref{B2}, noting that each 
configuration in $B_{\st,L}$ is monochromatic, so that 
$|B_{\st,L}|=2\,q$.
This shows that also $p_\st^{z,q}\le q^{-1/12}$. Finally, for $B_0$ 
we  compare the set $B_{0,L}$ with 
\be{F1L}
F_{1,L}=\{\s\in\O_{L,\per}: \s_i\ne 0\mbox{ iff }i\in 2\,\Z^2\}\;;
\ee
this gives $p_0^{z,q}\le 
(zq)^{-1/4}$. The counterpart of Proposition \ref{contour-est} thus holds 
as soon as $q_0$ is so large that $q_0^{-1/56}+ 
2\,q_0^{-1/12}+q_0^{-1/4}\le\d$. 

With the contour estimate in hand we can then proceed as in Section 
\ref{sec:compete}. Proposition \ref{phases} carries over verbatim; 
the only difference is that $G_\st$ is replaced by $G_\dis$ (which is 
not divided into two parts with disjoint side-projections), and 
$\r(\mu_\dis)\le \frac 14(1-2\e)+2\e=\frac 14 +\frac{3\e}{2}$. 
By the latter estimate, the
assumption $\e<1/6$ is slightly stronger than necessary
for adapting Lemma \ref{z_c} to the present case,
but we stick to it for simplicity.

The counterparts of Lemmas \ref{stag} and \ref{order} are obtained as
follows.  On the one hand, we have the estimate
\[
\mu_{L,\per}^{z,q}(G_{\ord,L})\le 
\mu_{L,\per}^{z,q}(G_{\ord,L})\Big/\mu_{L,\per}^{z,q}(F_{1,L})
\le z^{v(L)}q \;z^{-v(L)/4}q^{-v(L)/4}\;,
\]
showing that 
\[
\mu\Big(0\in V(G_\ord)\Big)\le (z^3/q)^{1/4}
\le 3^{-3/4}< (1-\d)/2
\]
when $\mu\in\bar\G_\per(z,q)$ and $z\le q^{1/3}/3$.  On the other hand,
as in Lemma \ref{B10} we find
\[
\mu_{L,\per}^{z,q}(G_{\dis,L})\le 
\mu_{L,\per}^{z,q}(G_{\dis,L})\Big/\mu_{L,\per}^{z,q}(G_{\ord,L})
\le \frac{2^{14L+2} z^{v(L)/4}q^{v(L)/4}}{z^{v(L)}q}
\]
and therefore
\[
\mu\Big(0\in V(G_\dis)\Big)\le (q/z^3)^{1/4}
\le 3^{-3/4}< (1-\d)/2
\]
when $\mu\in\bar\G_\per(z,q)$ and $z\ge 3\,q^{1/3}$.  With these
ingredients it is now straightforward to complete the proof of 
Theorem \ref{th:sqWR} along the lines of Section
\ref{sec:compete} .

\section{Proof of Theorem \ref{th:WRhc}}
\label{sec:pf_hc}

Here we consider the Widom--Rowlinson model with molecular hard-core
exclusion.  We look again at good configurations in plaquettes.  The
set $\O_C$ of admissible configurations in $C$ splits into the
good sets
\[
G_\ord= G_\even\cup G_\odd =\bigcup_{1\le a\le q} G_{a,\even} \cup 
G_{a,\odd} \equiv \bigcup_{1\le a\le q}\{\textstyle{0\;a \choose a\;0}\}
\cup \{\textstyle{a\;0 \choose 0\;a}\}
\]
of \emph{ordered staggered} configurations, the good set
\[ 
G_\dis =B_1=\{{\textstyle{0\;0 \choose a\;0}, {0\;0 \choose 0\;a}, 
{0\;a \choose 0\;0},{a\;0 \choose 0\;0}}: 1\le a\le q\}
\]
of \emph{disordered} configurations, and the only bad set $B_0$ 
consisting of the empty configuration. The main technical problem 
which is new in the present model is that the sets $G_\ord$ and 
$G_\dis$ \emph{fail} to have disjoint side-projections (although this 
is the case for the sets $G_{a,\even}$ and $G_{a,\odd}$). We 
therefore cannot simply consider sets of good plaquettes, but need to 
consider the sets of ``good plaquettes with neighbors in the same phase''.
Accordingly, we introduce the random sets
\bea 
\hat V(G_\ord)&=& \{i\in V(G_\ord):i+(1,0),i+(0,1)\in V(G_\ord)\}\;,\\
\hat V(G_\dis)&=& \{i\in V(G_\dis):i+(1,0),i+(0,1)\in V(G_\dis)\}\;,
\eea
and $\hat V(G)=\hat V(G_\ord)\cup\hat V(G_\dis)$.  By definition, a
sea in $\hat V(G)$ then contains either a sea in $\hat V(G_\ord)$ or a
sea in $\hat V(G_\dis)$.  To establish the existence of such a sea we
use the following contour estimate.
%
\begin{prop}\label{contour-est-hc} For any $\d>0$ there exists a 
number $q_0\in\N$ such that 
\[
\mu(\D \cap \hat V(G)=\emptyset) \le \d^{|\D |}
\]
whenever $q\ge q_0$, $zq\ge q_0^{1/7}$, $\mu\in \G_\per(z,q)$, and 
$\D \subset \Z^2$ is finite.
\end{prop}
\proof Let us start by introducing some notations. We consider the 
sublattices
\[
\LL_{1,\even}=\{i=(i_1,i_2)\in \Z^2: i_1\mbox{ is even}\}\;,\quad 
\LL_{1,\odd}= \Z^2\setminus \LL_{1,\even}\;,
\]
and their rotation images $\LL_{2,\even}$ and $\LL_{2,\odd}$ which 
are similarly defined. We also introduce the horizontal 
double-plaquette
\[ 
D_1=C\cup (C+(1,0)) =\{0,1,2\}\times\{0,1\}
\]
and the event
\[ 
E_1 =\Big\{ \s\in E^{D_1}: \s_C\in G_\ord,\, \s_{C+(1,0)}\in G_\dis,
\mbox{ or vice versa}\Big\}
\]
that the two sub-plaquettes of $D_1$ are good but of different type.
$E_1$ thus consists of the configurations of the form ${a\;0\;0 
\choose 0\;a\;0}$ with $1\le a\le q$, and their reflection images.
In the same way, we define the vertical double-plaquette $D_2=C\cup 
(C+(0,1))$ and the associated event $E_2$. With these notations we 
have
\[ 
\Z^2\setminus \hat V(G) \subset \bigcup_{k=1}^7 W_k\;,
\]
where the random subsets $W_k$ of $\Z^2$ are given by
\bea
W_1 &=& V(B_0)\;,\quad  W_2\ =\ V(B_0)-(1,0)\,,\quad  W_3\ =\ V(B_0)-(0,1)\,,\\
W_4 &=& \{i\in \LL_{1,\even}: \s_{D_1+i}\in E_1 \}\;,\quad
W_5\ =\ \{i\in \LL_{1,\odd} : \s_{D_1+i}\in E_1 \}\;,\\
W_6 &=& \{i\in \LL_{2,\even}: \s_{D_2+i}\in E_2 \}\;,\quad
W_7\ =\ \{i\in \LL_{2,\odd} : \s_{D_2+i}\in E_2 \}\;.
\eea
(The sets $W_k$ are not necessarily disjoint.) So, for each 
$\mu\in\G_\per(z,q)$ we can write
\be{cases}
\mu(\D\cap\hat V(G)=\emptyset) \le \sum_{\D_1\cup\ldots\cup 
\D_7=\D } 
\min_{1\le k\le 7}\mu(\D_k\subset W_k)\;,
\ee
where the sum extends over all disjoint partitions of $\D$. 
We estimate now each term.

Consider first the case $k=1$. Just as in Lemma \ref{B10} we obtain 
from the chessboard estimate
\[ 
\mu(\D_1\subset W_1)^{1/|\D_1|}\le \limsup_{L\ti} 
\frac{\mu_{L,\per}^{z,q}(B_{0,L})^{1/v(L)}}{\mu_{L,\per}^{z,q}
(F_{1,L})^{1/v(L)}}= (zq)^{-1/4}\;,
\]
where $F_{1,L}$ is given by \rf{F1L}. The same estimate holds in the 
cases $k=2,3$ because these merely correspond to a translation.

Next we turn to the case $k=4$. Let $L$ be so large that $\L_L\supset 
\D_4$. Using reflection positivity in the lines through the sites 
of $\LL_{1,\even}$, we conclude from the chessboard estimate that
\[ 
\mu_{L,\per}^{z,q}(\D_4\subset W_4)^{1/|\D_4|}\le  
\mu_{L,\per}^{z,q}(E_{1,L})^{2/v(L)}
\]
for the event
\[ 
E_{1,L} =\Big\{ \s\in \O_{L,\per}: \s_{D_1(i)}\in E_1
\mbox{ for all }i\in\L_L\cap\LL_{1,\even}\Big\}\;.
\]
In the above, $D_1(i)$  stands for the image $D_1+i \mbox{ mod }\L_L$ of 
$D_1$ under the periodic shift by $i$ of the torus $\L_L$.
Each $\s\in E_{1,L}$ has the following structure: every fourth 
vertical line (with horizontal coordinate either 0 or 2 
modulo 4) is empty, and on each group of three vertical lines 
between these empty lines every second site is occupied, with the 
coordinates of occupied sites being either even-odd-even in these 
three lines, or odd-even-odd; see the figure below.

\scriptsize
\[
\ba{*{9}{c}}
\c & \b & \c & \b & \c & \c & \b & \c & \c \\
\c & \c & \b & \c & \c & \b & \c & \b & \c \\
\c & \b & \c & \b & \c & \c & \b & \c & \c \\
\c & \c & \b & \c & \c & \b & \c & \b & \c \\
\c & \b & \c & \b & \c & \c & \b & \c & \c \\
\c & \c & \b & \c & \c & \b & \c & \b & \c \\
\c & \b & \c & \b & \c & \c & \b & \c & \c
\ea
\]\normalsize
Of course, the interaction implies that 
the color of particles is constant in each of these groups of three 
vertical lines. Consequently, each such $\s$ has particle number 
$3v(L)/8$, and $|E_{1,L}| = 2 (2q)^{6L}$; recall the definition 
\rf{LaL} of $\L_L$.

We now make a construction similar to that in Lemma \ref{B3}.  We
divide $\L_L$ into $12L^2$ rectangular cells $\D(j)$ of size
$8\times7$.  Let $\D_0(j)$ be the rectangular cell of size $7\times6$
situated in the left lower corner of $\D(j)$, and consider the set
\[
 F_L=\{\s\in\O_{L,\per}: \s_i\ne 0\mbox{ iff }i_1+i_2
\mbox{ is even and }i\in\D_0(j)\mbox{ for some }j\}\;.
\]
Since $|\D(j)\setminus\D_0(j)|=|\D(j)|/4$ for all $j$, each $\s\in
F_L$ has particle number $3\,v(L)/8$, just as the configurations in
$E_{1,L}$.  As the colors of the particles in the blocks $\D_0(j)$ can
be chosen independently, we have $|F_L|=q^{12L^2}=q^{v(L)/56}$. Hence
\[
\mu_{L,\per}^{z,q}(E_{1,L}) \le 
\frac{\mu_{L,\per}^{z,q}(E_{1,L})}{\mu_{L,\per}^{z,q}(F_L)}
= \frac{|E_{1,L}|}{|F_L|} \le 2 (2q)^{6L}\; q^{-v(L)/56}
\]
and therefore, by taking the $2/v(L)$'th power and letting $L\ti$,
we obtain
\[ 
\mu(\D_4\subset W_4)^{1/|\D_4|}\le q^{-1/28}\;.
\]
The same estimate holds in the cases $k=5,6,7$, as these are obtained 
by a translation or interchange of coordinates.

We now combine all previous estimates as follows.  Let $q_0$ be so
large that $7\, (q_0^{-1/28})^{1/7}<\d$, and suppose that $q\ge q_0$
and $zq\ge q_0^{1/7}$.  Then $\mu(\D_k\subset W_k)\le
q_0^{-|\D_k|/28}$ for all $k$ and thus, in view of \rf{cases} and
since $|\D_k|\ge |\D|/7$ for at least one $k$,
\[
\mu(\D\cap\hat V(G)=\emptyset) \le 
\sum_{\D_1\cup\ldots\cup \D_7=\D } (q_0^{-1/28})^{|\D|/7}
<\d^{|\D|}.
\]
The proof of the contour estimate is therefore complete.  
$\Box$

\medskip
To prove Theorem \ref{th:WRhc} we can now proceed as in Section 
\ref{sec:compete}. Let $\hat\sea(G)$ be the largest sea in $\hat V(G)$ 
if the latter contains a sea, and $\hat\sea(G)=\emptyset$ otherwise. 
It is then immediate that a counterpart of Proposition \ref{sea} 
holds, and the definition of $\hat V(G)$ implies that
\[
\{\hat\sea(G)\ne\emptyset\}=\{\hat\sea(G_\dis)\ne\emptyset\}\cup 
\{\hat\sea(G_\ord)\ne\emptyset\}\;,
\]
where the two sets on the right-hand side are disjoint. Moreover,
\[ 
\{\hat\sea(G_\ord)\ne\emptyset\}\subset
\bigcup_{a=1}^q\{\sea(G_{a,\even})\ne\emptyset\}
\cup \{\sea(G_{a,\odd})\ne\emptyset\}\;.
\]
By the argument of Proposition \ref{phases} we thus obtain the 
existence of $2q$ ordered phases (as described in Theorem 
\ref{th:WRhc}(i)) whenever $z$ is such that 
$\mu(0\in V(G_\ord))\ge \mu(0\in V(G_\dis))$
for some $\mu\in \bar\G_\per(z,q)$, and the existence of a disordered 
phase $\mu_\dis$ whenever the reverse inequality holds for such a $\mu$.
We have $\r(\mu_\dis)\le \frac 14 + \frac{3\e}{2}$ and 
$\r(\mu_{a,\even})= \r(\mu_{a,\odd})\ge \frac 12 -\e$. 
The topological argument of Section \ref{sec:compete} together with 
obvious counterparts of Lemmas \ref{stag} and \ref{order} then show
that both cases must occur simultaneously for some $z=z_c(q)$, and 
this $z$ is unique by the convexity argument of Lemma \ref{z_c}. 
(For the latter we need to assume that $\e<1/10$.)

We conclude this section with a comment on the model with 
nearest-particle color repulsion and a molecular hard-core exclusion
between next-nearest neighbors. 

\medskip\noindent
{\bf Comment on Remark \ref{th:WRhc} (3).}
If the r\^oles of $\Phi$ and $U$ are interchanged, the good ordered 
configurations in 
$C$ are those with two particles of the same color on one side 
of $C$, and no particle on the opposite side; we call this set again 
$G_\ord$. For large $z$, one can 
easily establish a contour estimate implying the existence of a sea 
$S(G_\ord)$, and thus by symmetry also the existence of the 
four phases mentioned in Remark \ref{th:WRhc} (3). The disordered good 
plaquettes are again described by the set $G_\dis$. 
As in the case of the Hamiltonian \rf{Ham_hc}, the sets $G_\ord$ and 
$G_\dis$ have no disjoint side-projections. However, whereas in 
that case we were able to show an entropic disadvantage 
in having adjacent $G_\ord$- and 
$G_\dis$-plaquettes, this is not true in the present case. The 
configurations resulting from iterated reflections of a double 
plaquette of type ${\c\,\c\,\c\choose\b\,\b\,\c}$ have the maximal 
entropy possible for this particle number. Therefore the system 
can freely combine ordered and disordered 
plaquettes, and our argument for a first-order transition breaks down.
So it seems likely that the transition from the ordered to the 
disordered phase is of second order.

\section{Proof of Theorem \ref{th:contspin}}
\label{sec:pf_cs}

The analysis of the plane-rotor Widom--Rowlinson model is very similar
to that of the standard Widom--Rowlinson model; only a few
modifications are necessary.  We define again the set $\O_C$ of
admissible configurations in the plaquette $C$ in the obvious way,
introduce the sets $G_\st=G_\even\cup G_\odd$ as in Section
\ref{sec:compete} (replacing $\{1,\ldots,q\}$ by $S^1$), and set
$G_\ord=\{\s\in \O_C: \s_i\in S^1 \mbox{ for all }i\in C\}$ and
$G=G_\ord\cup G_\st$.  The main task is to obtain a counterpart of the
contour estimate, Proposition \ref{contour-est}.  To this end we
consider the same classes $B_k$, $k\in\{0,\ldots,3\}$ of bad
configurations as in Section \ref{sec:cont_est} (with the obvious
modifications), and the sets $B_{k,L}$ and the associated quantities
$p_k^{z,\a}$.

To deal with the case $k=3$ we proceed as in Lemma \ref{B3}, arriving 
at the inequality
\[
\mu_{L,\per}^{z,\a}(B_{3,L}) \le 
\frac{\mu_{L,\per}^{z,\a}(B_{3,L})}{\mu_{L,\per}^{z,\a}(F_{3,L})}
= \frac{\nu^\L(B_{3,L})}{\nu^\L(F_{3,L})}\;.
\]
Now, $\nu^\L(B_{3,L})\le 2^{14L+2} (2\a)^{3\,v(L)/4-1}$; 
the first factor estimates the number of possible occupation 
patterns, and the second term bounds the probability that the 
configuration is admissible (by keeping only the bonds in a tree 
spanning all occupied positions). On the other hand, 
$\nu^\L(F_{3,L})\ge (\a^{7\cdot 6-1})^{v(L)/56}$, as can be seen 
by letting the spins in each block $\D_0(k)$ follow a ``leader spin'' up 
to the angle $2\pi\a/2$. Hence
\[
\frac{\nu^\L(B_{3,L})}{\nu^\L(F_{3,L})}\le 
2^{14L+2}\;2^{3\,v(L)/4}\;\a^{v(L)/56-1}
\]
and therefore $p_3^{z,\a}\le 2^{3/4}\,\a^{1/56}$. 

In the case $k=2$ we proceed as in the proof of Lemma \ref{B2}. On the one 
hand, 
\[
\nu^\L(B_{2,L})\le 
2((2\a)^{24\,L-1})^{14\,L}+2((2\a)^{28\,L-1})^{12\,L}
\le 4 (2\a)^{v(L)/2-14L}
\]
since the spins are ordered in separate rows or columns, and $2\a<1$. On the 
other hand, $\nu^\L(F_{2,L})\ge (\a^{2\cdot 3-1})^{v(L)/12}$ by the 
same argument as above. Hence
\[
\mu_{L,\per}^{z,\a}(B_{2,L}) \le 
\frac{\nu^\L(B_{2,L})}{\nu^\L(F_{2,L})}
\le 4\cdot 2^{v(L)/2}\; \a^{v(L)/12}\; (2\a)^{-14L}
\]
and therefore $p_2^{z,\a}\le 2^{1/2}\,\a^{1/12}$.

Finally, for $k=0$ we obtain 
\[
\mu_{L,\per}^{z,\a}(B_{0,L}) \le 
\frac{\mu_{L,\per}^{z,\a}(B_{0,L})}{\mu_{L,\per}^{z,\a}(G_{\even,L})}
= \frac{1}{z^{v(L)/2}}
\]
and thus $p_0^{z,\a}\le z^{-1/2}$. Likewise, in the case $k=1$ we get
as in Lemma \ref{B10}
\[
\mu_{L,\per}^{z,\a}(B_{1,L}) \le 
\frac{\mu_{L,\per}^{z,\a}(B_{1,L})}{\mu_{L,\per}^{z,\a}(G_{\even,L})}
\le \frac{2^{14L+2}z^{v(L)/4}}{z^{v(L)/2}}
\]
and thereby $p_1^{z,\a}\le z^{-1/4}$. Combining these estimates as in the 
proof of Proposition \ref{contour-est} we arrive at the counterpart 
of \rf{cont-est} as soon as $\a_0$ is so small that
\[
2^{3/4}\,\a_0^{1/56}+2^{1/2}\,\a_0^{1/12}+\a_0^{1/4}+\a_0^{1/2}\le\d
\]
and $\a\le\a_0$, $z\ge 1/\a_0$.  

To complete the proof of Theorem \ref{th:contspin} as in Section
\ref{sec:compete} we still need to adapt Lemmas \ref{stag} and
\ref{order}.  Writing
\[
\mu_{L,\per}^{z,\a}(G_{\ord,L})\le 
\frac{\mu_{L,\per}^{z,\a}(G_{\ord,L})}{\mu_{L,\per}^{z,\a}(G_{\even,L})}
\le \frac{z^{v(L)}(2\a)^{v(L)-1} }{z^{v(L)/2}}
\]
we find that for $z\le \a^{-2}/18$ and $\mu\in\bar\G_\per(z,\a)$
\[
\mu\Big(0\in V(G_\ord)\Big)\le z^{1/2}2\a\le 2/\sqrt{18} < 
(1-\d)/2\;.
\]
Likewise, since
\[
\mu_{L,\per}^{z,\a}(G_{\st,L})\le 
\frac{\mu_{L,\per}^{z,\a}(G_{\st,L})}{\mu_{L,\per}^{z,\a}(G_{\ord,L})}
\le \frac{2\;z^{v(L)/2} }{z^{v(L)}\,\a^{v(L)-1}}\;,
\]
we see that for $z\ge 5\,\a^{-2}$ and $\mu\in\bar\G_\per(z,\a)$
\[
\mu\Big(0\in V(G_\st)\Big)\le z^{-1/2}\a\inv < 
(1-\d)/2\;.
\]
The remaining arguments of Section \ref{sec:compete} can be taken over
with no change to prove Theorem \ref{th:contspin}.


\end{document}